\documentclass[12pt]{article}
\usepackage{mathrsfs}
\usepackage[hypertex]{hyperref}
\usepackage{amsfonts}
\usepackage{graphicx}
\usepackage{mathptmx}
\usepackage{latexsym,amsmath,amssymb,amsfonts,amsthm}

\newtheorem{theorem}{Theorem}[section]
\newtheorem{lemma}[theorem]{Lemma}

\newtheorem{corollary}[theorem]{Corollary}

\newtheorem{remark}[theorem]{Remark}

\newtheorem{defn}[theorem]{Definition}

\begin{document}
\setcounter{page}{1}
\title{Open manifold with nonnegative Ricci curvature and collapsing volume}
\author{Jing Mao}
\date{}
\protect\footnotetext{\!\!\!\!\!\!\!\!\!\!\!\!{MSC 2010: 53C20;
53C21}
\\
{ ~~Key Words: Critical point; Collapsing volume; Radial sectional curvature; Finite topological type.} \\
 supported by Funda\c{c}\~{a}o de Ci\^{e}ncia e
Tecnologia (FCT) through a doctoral fellowship SFRH/BD/60313/2009.}
\maketitle ~~~\\[-15mm]
\begin{center}{\footnotesize  Departamento de Matem\'{a}tica,
Instituto Superior T\'{e}cnico, Technical University of Lisbon,
Edif\'{\i}cio Ci\^{e}ncia, Piso 3, Av.\ Rovisco Pais, 1049-001
Lisboa, Portugal; jiner120@tom.com}
\end{center}

\begin{abstract} In this paper, an $n$-dimensional complete open manifold
with nonnegative Ricci curvature and collapsing volume has been
investigated. If its radial sectional curvature bounded from below,
it shows that such a manifold is of finite topological type under
some restrictions shown below.
\end{abstract}

\markright{\sl\hfill  J. Mao \hfill}

\section{Introduction}
\renewcommand{\thesection}{\arabic{section}}
\renewcommand{\theequation}{\thesection.\arabic{equation}}
\setcounter{equation}{0} \setcounter{maintheorem}{0}

Without specification, in this paper we let $(M,g)$ be an
$n$-dimensional complete Riemannian manifold with nonnegative Ricci
curvature, denote by $B(p,r)$ the open geodesic ball centered at a
 point $p\in{M}$ with radius $r$ and ${\rm{vol}}(B(p,r))$ its volume,
and let $w_{n}$ be the volume of unit ball in the Euclidean space
$R^{n}$. By the classical volume comparison theorem \cite{rr,mjp},
we know the function
$r\rightarrow\frac{{\rm{vol}}(B(p,r))}{w_{n}r^{n}}$ is monotone
decreasing. Define
\begin{eqnarray*}
\alpha_{M}:=\lim\limits_{r\rightarrow\infty}\frac{{\rm{vol}}(B(p,r))}{w_{n}r^{n}},
\end{eqnarray*}
it is not difficult to prove $\alpha_{M}$ is independent of the
choice of $p$, which implies $\alpha_{M}$ is a global geometric
invariant. Obviously, $\alpha_{M}\in[0,1]$, and
\begin{eqnarray} \label{1.1}
\alpha_{M}w_{n}r^{n}\leq{{\rm{vol}}(B(p,r))}\leq{w_{n}r^{n}}, \quad
for ~\forall{p}\in{M}~~and~~\forall{r}>0.
\end{eqnarray}

We say that $(M,g)$ has large volume growth provided $\alpha_{M}>0$.
Riemannian manifold with nonnegative Ricci curvature and large
volume growth has been investigated intensively and some good
results have been obtained in the past decades. Let $(N,g)$ be an
$n$-dimensional complete open manifold with Ricci curvature
$Ric_{N}\geq0$ and $\alpha_{N}>0$. By Bishop-Gromov comparison
theorem \cite{rr,mjp}, $N^{n}$ is isometric to $R^{n}$ when
$\alpha_{N}=1$. It has been shown by Li \cite{p} that $N$ has finite
fundamental group. Anderson \cite{ma} has proved that the order of
$\pi_{1}(N)$ is bounded from above by $\frac{1}{\alpha_{N}}$.
Petersen has conjectured that if $\alpha_{N}>\frac{1}{2}$, then
$N^{n}$ is diffeomorphic to $R^{n}$ in \cite{pe}. A theorem has been
proved by do Carmo and Xia in \cite{dx} to support this conjecture.
Xia \cite{cx1,cx2} has shown that $N^{n}$ is diffeomorphic to a
Euclidean space $R^{n}$ or is of finite topological type under
different restrictions about $\alpha_{N}$ and volume ratio
$\frac{{\rm{vol}}(B(p,r))}{w_{n}r^{n}}$.

It seems that if we want to get conclusions like $M^{n}$ is
diffeomorphic to a Euclidean space $R^{n}$ or is of finite
topological type for an $n$-dimensional complete open Riemannian
manifold $(M,g)$ with nonnegative Ricci curvature, the large volume
growth condition can not be avoided. However, in this paper, we find
that we could use the $\alpha$-order collapsing volume condition to
replace the large volume growth condition. Actually, for $M^{n}$
with nonnegative Ricci curvature, we have
\begin{eqnarray*}
c_{1}(n)vol(B(p,1))r\leq{vol(B(p,r))}\leq{w_{n}r^{n}},
\end{eqnarray*}
where $c_{1}(n)$ is a constant depending on $n$. This inspires us
that maybe we could consider some restriction on
${\rm{vol}}(B(p,1))$ if we want to get conclusions similar with
those in \cite{cx1}, since ${\rm{vol}}(B(p,1))$ has connection with
the volume ratio $\frac{{\rm{vol}}(B(p,r))}{w_{n}r^{n}}$.
Fortunately, in this paper we prove that this is possible. Now, we
want to give our main result in this paper, however, before that we
need to introduce some concepts and conclusions. First, we need to
use the following notion introduced in \cite{z,zg}

\begin{defn} \label{def1} Let $M$ be a complete noncompact manifold
and let $p\in{M}$ be a point such that
\begin{eqnarray} \label{1.2}
v_{p}(r)=\inf\limits_{x\in{S(p,r)}}{\rm{vol}}(B(x,1))=\emph{O}\left(\frac{1}{r^{\alpha}}\right),
\end{eqnarray}
 where $S(p,r)$ denotes the geodesic sphere centered at $p$ with
radius $r$ on $M$, then we say that $M$ has $\alpha$-order
collapsing volume.
\end{defn}

We also need the following lemma in \cite{zg}
\begin{lemma} \label{lemma1}
Let $M$ be a complete noncompact n-manifold with nonnegative Ricci
curvature $Ric_{M}\geq0$. Then there is a constant $c_{2}$ such that
for $\forall{R}\geq{r}$, we have
\begin{eqnarray} \label{1.3}
vol(B(p,R)\setminus{B(p,r)})\leq{c_{2}}\int_{r}^{R}\frac{{\rm{vol}}(B(p,s))}{s}ds.
\end{eqnarray}
\end{lemma}

A manifold M is said to have finite topological type if there is a
compact domain $\Omega$ whose boundary $\partial\Omega$ is a
topological manifold such that $M\setminus\Omega$ is homeomorphic to
$\partial\Omega\times[0,\infty)$. For a fixed point $p\in{M}$, we
say its radial sectional curvature, $K_{p}^{min}$, bounded from
below by a constant $-C$ if for any minimal geodesic $\gamma$
starting from $p$ all sectional curvatures of the planes which are
tangent to $\gamma$ are greater than or equal to $-C$, i.e.
$K_{p}^{min}\geq-C$. The main result is the following

\begin{theorem} \label{theorem2} Let $(M,g)$ be an $n$-dimensional
complete noncompact Riemannian manifold with nonnegative Ricci
curvature $Ric_{M}\geq0$ and $\alpha$-order collapsing volume
($0\leq\alpha\leq\frac{1}{n}$). Suppose that $K_{p}^{min}\geq-C$ for
some point $p\in{M}$ and some positive constant $C$. If
\begin{eqnarray} \label{1.5}
\limsup\limits_{r\rightarrow\infty}\left[\frac{{\rm{vol}}(B(p,r))}{r^{1+\frac{1}{n}-\alpha}}\right]
\leq{c_{3}}\left(\frac{\log2}{8\sqrt{C}}\right)^{\frac{n-1}{n}},
\end{eqnarray}
where $c_{3}$ is a positive constant depending on $c_{2}$, then $M$
has finite topological type.
\end{theorem}

 Our Theorem \ref{theorem2} is a generalization of Zhan's main theorem 3 in \cite{z}, since
 we only need the radial sectional curvature $K_{p}^{min}$ bounded from below. Besides, we
 will give a more general version of our main theorem in the last
 section, which is a generalization of theorem 8 in \cite{z} and shows the advantage of our result indeed.
 The paper is organized as follows. Some useful conclusions will be
introduced and proved in Section 2. In Section 3, we will give the
proof of Theorem \ref{theorem2}.

\section{Some useful facts}
\renewcommand{\thesection}{\arabic{section}}
\renewcommand{\theequation}{\thesection.\arabic{equation}}
\setcounter{equation}{0} \setcounter{maintheorem}{0}

First, we would like to give an isotopy lemma obtained by Grove and
Shiohama which will be used in the proof of our main theorem later.

\begin{lemma} (\cite{kk}) \label{lemma2} If
$r_{1}\leq{r_{2}}\leq\infty$ and a connected component $D$ of
$\overline{B(p,r_{2})}\setminus{B(p,r_{1})}$ is free of critical
points of $p$, then $D$ is homeomorphic to
$D_{1}\times[r_{1},r_{2}]$, where $D_{1}$ is a topological
submanifold without boundary.
\end{lemma}

For convenience, throughout this paper, all geodesics are assumed to
have unit speed. In order to have a topological cognizance of the
above lemma, we want to recall the notion of critical point here.
For a point $p\in{M}$, let $d_{p}(x)=d(p,x)$, where $d$ is the
metric on the Riemannian manifold $M$, obviously, the function
$d_{p}$ is Lipschitz continuous, however, it is not a smooth
function on the cut locus of $p$, which implies the critical points
of $d_{p}$ can not be defined in a usual way. The notion of critical
points of $d_{p}$ was introduced minutely in \cite{kk}. A point
$q\in{M}$ different from $p$ is  called a critical point of $d_{p}$
if there always exists a minimizing geodesic $\gamma$ from $q$ to
$p$ such that for any $v\in{T_{q}M}$, the forming angle
$\angle(v,\gamma'(0))$ satisfies
$\angle(v,\gamma'(0))\leq\frac{\pi}{2}$. We simply say $q$ is a
critical point of $p$. By the above isotopy lemma, we know that an
$n$-dimensional complete noncompact Riemannian manifold $M$ is
diffeomorphic to a Euclidean space $R^{n}$ if there is a point
$p\in{M}$ such that $p$ has no critical points other than $p$, which
shows the importance of this lemma.

Now, we want to recall a notion named $k$-th Ricci curvature
($1\leq{k}\leq{n-1}$) for the $n$-dimensional Riemannian manifold
$M^{n}$. We say that the $k$-th Ricci curvature of $M$ is
nonnegative if for any point $x\in{M}$ and any mutually orthogonal
unit tangent vector $e,e_{1},\cdots,e_{k}\in{T_{x}M}$, we have
$\sum_{i=1}^{k}K(e\wedge{e_{i}})\geq0$, here $K(e\wedge{e_{i}})$ is
the sectional curvature of the plane spanned by $e$ and $e_{i}$.
Denote this fact by $Ric_{M}^{(k)}\geq0$. Notice that if
$Ric_{M}^{(k)}\geq0$, then $Ric_{M}\geq0$. Let $p,q\in{M}$, then the
excess function is defined by
\begin{eqnarray*}
e_{pq}(x):=d(p,x)+d(q,x)-d(p,q).
\end{eqnarray*}
We have the following lemma which gives an upper bound for the
excess function.

\begin{lemma} (\cite{kd,zs}) \label{lemma3} Let $(M,g)$ be an $n$-dimensional
complete Riemannian manifold with $Ric_{M}^{(k)}\geq0$ for some
$1\leq{k}\leq{n-1}$. Let $\gamma:[0,a]\rightarrow{M}$ be a minimal
geodesic from $p$ to $q$. Then for any $x\in{M}$,
\begin{eqnarray} \label{2.1}
e_{pq}(x)\leq8\left(\frac{s^{k+1}}{r}\right)^{\frac{1}{k}},
\end{eqnarray}
where $s=d(x,\gamma)$, $r=\min(d(p,x),d(q,x))$.
\end{lemma}

Let $\gamma:[0,\infty)\rightarrow{M}$ be a ray emanating from $p$.
For any $x\in{M}$, by triangle inequality, it is easy to see that
$e_{p,\gamma(t)}(x)=d(p, x)+d(\gamma(t),x)-t$ is decreasing in $t$
and that $e_{p,\gamma(t)}(x)\geq0$. Define the excess function
$e_{p,\gamma}$ associated to $p$ and $\gamma$ as
\begin{eqnarray} \label{2.2}
e_{p,\gamma}(x)=\lim\limits_{t\rightarrow\infty}e_{p,\gamma(t)}(x).
\end{eqnarray}
Obviously, $e_{p,\gamma}(x)\leq{e}_{p,\gamma(t)}(x)$ for any $t>0$.
For this function $e_{p,\gamma}$, Xia \cite{cx1} has proved the
following lemma.

\begin{lemma}  \label{lemma4} Let $(M,g)$ be a complete open
Riemannian manifold with $K_{p}^{min}\geq-C$ for some $C>0$ and
$p\in{M}$. Suppose $x\neq{p}$ is a critical point of $p$. Then for
any ray $\gamma:[0,\infty)\rightarrow{M}$ issuing from $p$, we have
\begin{eqnarray} \label{2.3}
e_{p,\gamma}(x)\geq\frac{1}{\sqrt{C}}\log\left(\frac{2}{1+e^{-2\sqrt{C}d(p,x)}}\right).
\end{eqnarray}
\end{lemma}

At the end of this section, we want to give a lemma which will play
an important role in the proof of our main theorem. In order to
prove the lemma, we have to give some notions first. Let $M$ be an
$n$-dimensional complete open Riemannian manifold. For a given point
$p\in{M}$, set
\begin{eqnarray} \label{2.4}
v_{p}(A,r)=\inf\limits_{x\in{S(p,r)}}{\rm{vol}}\left(B\left(x,\frac{A}{2}\right)\right),
\quad\quad 0<A<\frac{r}{2},
\end{eqnarray}
here $S(p,r)$ has the same meaning as before. For $r>0$ and a point
$p\in{M}$, let
 \begin{eqnarray*}
 R(p,r)=\{\gamma(r)|\gamma~\emph{is a ray from p}\},
 \end{eqnarray*}
 obviously, $R(p,r)$ is the set of points of the intersections of
 the geodesic sphere centered at $p$ of radius $r$ with all the rays
 issuing from $p$. Let
 \begin{eqnarray} \label{2.5}
R_{p}(x)=d(x,R(p,r)), \quad {\rm{where}}~~r=d(p,x).
 \end{eqnarray}
  We can prove
\begin{lemma} \label{lemma5} Let $M$ be an $n$-dimensional complete
noncompact Riemannian manifold with nonnegative Ricci curvature
$Ric_{M}\geq0$, and let $p\in{M}$, then for any $r>4$ and
$x\in{S(p,r)}$, we have
\begin{eqnarray} \label{2.6}
d(x,R_{p})\leq{c_{4}}\frac{{\rm{vol}}(B(p,r))\cdot(r+2)^{n}}{v_{p}(r)\cdot{r}^{n}}\log\left(\frac{r+2}{r-2}\right),
\end{eqnarray}
where $c_{4}=8c_{2}$ is a positive constant depending only on
$c_{2}$, $v_{p}(r):=v_{p}(2,r)$, and $R_{p}$ denotes the union of
rays issuing from $p$.
\end{lemma}

\begin{proof} Here we use a similar method as that of lemma 6 in \cite{z}. Let $\Omega_{r}$ be a boundary component of
$M\setminus{\overline{B(p,r)}}$ with
$\Omega_{r}\cap{R}(p,r)\neq\emptyset$. So, there exists a ray
$\gamma_{p}$ such that $\gamma_{p}(r)\in\Omega_{r}$. Let
$\{B(p_{j},\frac{A}{2})\}$ be a maximal set of disjoint balls with
radius $\frac{A}{2}$ and center $p_{j}\in\Omega_{r}$, here
$0<A<\frac{r}{2}$, then
\begin{eqnarray*}
\bigcup\limits_{j=1}^{N}B(p_{j},A)\supset\Omega_{r}
\end{eqnarray*}
and
\begin{eqnarray*}
N\leq\frac{{\rm{vol}}(B(p,r+A)\setminus{B}(p,r-A))}{v_{p}(A,r)},
\end{eqnarray*}
where $v_{p}(A,r)$ is defined as (\ref{2.4}). By the connectedness
of $\Omega_{r}$, we know that for any point $x\in\Omega_{r}$, there
exists a subset of $\{p_{j}\}_{j=1,\cdots,N}$, say
$\{q_{1},\cdots,q_{k}\}$, $k\leq{N}$, such that
\begin{eqnarray*}
B(q_{l},A)\cap{B}(q_{l+1},A)\neq\emptyset, \quad 1\leq{l}\leq{k-1},
\end{eqnarray*}
and $x\in{B(q_{1},A)}$, $\gamma_{p}(r)\in{B(q_{k},A)}$. Hence, we
can easily construct a piecewise smooth geodesic $c$ joining $x$ and
$\gamma_{p}(r)$ through $q_{l}$'s. This implies
\begin{eqnarray} \label{2.7}
d(x,\gamma_{p}(r))\leq{L(c)}\leq4NA\leq4\cdot\frac{{\rm{vol}}(B(p,r+A)\setminus{B}(p,r-A))}{v_{p}(A,r)}\cdot{A},
\end{eqnarray}
where $L(c)$ denotes the length of $c$. Then, by (\ref{2.5}), we
have
\begin{eqnarray*}
R_{p}(x)\leq4\cdot\frac{{\rm{vol}}(B(p,r+A)\setminus{B}(p,r-A))}{v_{p}(A,r)}\cdot{A},
\end{eqnarray*}
moreover,
\begin{eqnarray} \label{2.8}
d(x,R_{p})\leq{R_{p}(x)}\leq4\cdot\frac{{\rm{vol}}(B(p,r+A)\setminus{B}(p,r-A))}{v_{p}(A,r)}\cdot{A},
\end{eqnarray}
since $R(p,x)$ is only a part of the point set $R_{p}$.

On the other hand, By Lemma \ref{lemma1}, we can obtain
\begin{eqnarray*}
\frac{{\rm{vol}}(B(p,r+A)\setminus{B}(p,r-A))}{v_{p}(A,r)}&\leq&{\frac{c_{2}}{v_{p}(A,r)}}\int^{r+A}_{r-A}\frac{{\rm{vol}}(B(p,s))}{s}ds\\
&\leq&c_{2}\cdot\frac{{\rm{vol}}(B(p,r+A))}{v_{p}(A,r)}\log\left(\frac{r+A}{r-A}\right),
\end{eqnarray*}
furthermore, together with (\ref{2.8}), we have
 \begin{eqnarray*}
d(x,R_{p})\leq{4c_{2}A}\cdot\frac{{\rm{vol}}(B(p,r+A))}{v_{p}(A,r)}\log\left(\frac{r+A}{r-A}\right).
 \end{eqnarray*}
 Choose $A=2$, then we get
 \begin{eqnarray} \label{2.9}
d(x,R_{p})\leq{c_{4}}\cdot\frac{{\rm{vol}}(B(p,r+2))}{v_{p}(r)}\log\left(\frac{r+2}{r-2}\right),
 \end{eqnarray}
 here $c_{4}=8c_{2}$. By the volume comparison theorem
(see \cite{rr,mjp}), the expression (\ref{2.9}) becomes
 \begin{eqnarray*}
d(x,R_{p})\leq{c_{4}}\cdot\frac{{\rm{vol}}(B(p,r+2))}{v_{p}(r)}\log\left(\frac{r+2}{r-2}\right)\leq
{c_{4}}\cdot\frac{{\rm{vol}}(B(p,r))\cdot(r+2)^{n}}{v_{p}(r)\cdot{r}^{n}}\log\left(\frac{r+2}{r-2}\right),
 \end{eqnarray*}
 which implies our lemma.
\end{proof}

\begin{remark} \rm{Here we would like to point out that Lemma \ref{lemma5} is still true if we reduce the condition
$Ric_{M}\geq0$ to $Ric_{M}^{min}\geq0$, the radial Ricci curvature
is nonnegative. This is because for a complete open Riemannian
manifold with nonnegative radial Ricci curvature, the function
$r\rightarrow\frac{{\rm{vol}}(B(p,r))}{w_{n}r^{n}}$ is monotone
non-increasing, which has been proved by Shiohama in \cite{k}.}
\end{remark}

\section{Proof of the main theorem}
\renewcommand{\thesection}{\arabic{section}}
\renewcommand{\theequation}{\thesection.\arabic{equation}}
\setcounter{equation}{0} \setcounter{maintheorem}{0}

In fact, we could prove the following more general theorem than
Theorem \ref{theorem2}.

\begin{theorem} Let $(M,g)$ be an $n$-dimensional
complete noncompact Riemannian manifold with $Ric_{M}^{(k)}\geq0$
($1\leq{k}\leq{n-1}$). Suppose that $K_{p}^{min}\geq-C$ for some
point $p\in{M}$ and some positive constant $C$. If
\begin{eqnarray} \label{3.1}
\limsup\limits_{r\rightarrow\infty}\left[\frac{{\rm{vol}}(B(p,r))}{r^{1+\frac{1}{k+1}}\cdot{v}_{p}(r)}\right]
\leq{c_{5}}\left(\frac{\log2}{8\sqrt{C}}\right)^{\frac{k}{k+1}},
\end{eqnarray}
where $c_{5}=2^{-5}c_{2}^{-1}$ is a positive constant depending only
on $c_{2}$, and $v_{p}(r):=v_{p}(2,r)$, then $M$ has finite
topological type.
\end{theorem}
\begin{proof} We use a similar method as that of theorem 2.2 in \cite{cx1} to prove
our theorem. By the isotopy Lemma \ref{lemma2}, we know that if we
want to prove the complete Riemannian manifold $M$ is of finite
topological type, we only need to show that there are no critical
points outside a compact subset with respect to a fixed point
$p\in{M}$. Take an arbitrary point $x(\neq{p})\in{M}$ and set
$r=d(p,x)$, which implies $x\in{S(p,r)}$. Since
$c_{5}=2^{-2}c_{4}^{-1}=2^{-5}c_{2}^{-1}$,
$\limsup\limits_{r\rightarrow\infty}\left(\frac{r}{r+2}\right)^{n-1}=1,$
\begin{eqnarray*}
\limsup\limits_{r\rightarrow\infty}(r+2)\log\left(\frac{r+2}{r-2}\right)=\lim\limits_{r\rightarrow\infty}\left(\frac{r+2}{r-2}\right)\cdot
\lim\limits_{r\rightarrow\infty}
\log\left(1+\frac{4}{r-2}\right)^{(r-2)}=4,
\end{eqnarray*}
then our assumption (\ref{3.1}) enables us to find a small number
$\epsilon>0$ and a sufficiently large $r_{1}>4$ such that for any
$r\geq{r_{1}}$, we have
\begin{eqnarray} \label{3.2}
\frac{{\rm{vol}}(B(p,r))}{r^{1+\frac{1}{k+1}}\cdot{v_{p}(r)}}<c_{4}^{-1}\left[\left(\frac{r}{r+2}\right)^{n-1}-\epsilon\right]
\left[\frac{1}{(r+2)\log(\frac{r+2}{r-2})}-\epsilon\right]\left(\frac{\log2}{8\sqrt{C}}-\epsilon\right)^{\frac{k}{k+1}}.
\end{eqnarray}
On the other hand, since
\begin{eqnarray*}
\lim\limits_{r\rightarrow\infty}\log\left(\frac{2}{1+e^{-2\sqrt{C}r}}\right)=\log2,
\end{eqnarray*}
there is a sufficiently large $r_{2}$ such that
\begin{eqnarray} \label{3.3}
\frac{\log\left(\frac{2}{1+e^{-2\sqrt{C}r}}\right)}{8\sqrt{C}}>\frac{\log2}{8\sqrt{C}}-\epsilon,
\quad\quad \forall{r}\geq{r_{2}}.
\end{eqnarray}
Let $r_{0}=\max\{r_{1},r_{2}\}$, then from (\ref{3.2}) and
(\ref{3.3}) we have
\begin{eqnarray}
\frac{{\rm{vol}}(B(p,r))}{r^{n}\cdot{v_{p}(r)}}<c_{4}^{-1}\frac{r^{\frac{1}{k+1}}}{(r+2)^{n}}
\left[\log\left(\frac{r+2}{r-2}\right)\right]^{-1}\cdot\left[\frac{1}{8\sqrt{C}}\log\left(\frac{2}{1+e^{-2\sqrt{C}r}}\right)\right]^{\frac{k}{k+1}},
\end{eqnarray}
 for any $r\geq{r_{0}}$. By Lemma \ref{lemma5}, we could obtain
 \begin{eqnarray}
d(x,R_{p})<{r^{\frac{1}{k+1}}}\cdot\left[\frac{1}{8\sqrt{C}}\log\left(\frac{2}{1+e^{-2\sqrt{C}r}}\right)\right]^{\frac{k}{k+1}},
 \end{eqnarray}
 for any $r\geq{r_{0}}$. So, we can find a ray
 $\gamma:[0,\infty)\rightarrow{M}$ emanating from $p$ and
 satisfying
 \begin{eqnarray} \label{3.6}
 s:=d(x,\gamma)<{r^{\frac{1}{k+1}}}\cdot\left[\frac{1}{8\sqrt{C}}\log\left(\frac{2}{1+e^{-2\sqrt{C}r}}\right)\right]^{\frac{k}{k+1}},
 \end{eqnarray}
for any $r\geq{r_{0}}$. We can find a point $q\in\gamma$ such that
$d(x,q)=d(x,\gamma)$, moreover, by (\ref{3.6}), $d(x,q)<r$.
Additionally, by triangle inequality, we know
\begin{eqnarray}
\min(d(p,x),d(\gamma(t),x))=r, \quad\quad \forall{t}\geq2r.
\end{eqnarray}
Therefore $q\in\gamma((0,2r))$ and $d(x,\gamma|_{[0,2r]})=s$. Then
by Lemma \ref{lemma3}, (\ref{3.6}), and the fact
$e_{p,\gamma}(x)\leq{e}_{p,\gamma(t)}(x)$ for any $t>0$, we can
obtain
\begin{eqnarray} \label{3.8}
e_{p,\gamma}(x)\leq{e_{p,\gamma(2r)}(x)}
\leq8\left(\frac{s^{k+1}}{r}\right)^{\frac{1}{k}}
<\frac{1}{\sqrt{C}}\log\left(\frac{2}{1+e^{-2\sqrt{C}r}}\right).
\end{eqnarray}
So, by Lemma \ref{lemma4} and (\ref{3.8}), $x$ is not a critical
point of $p$. This implies $p$ has no critical point out of a
compact subset $\overline{B(p,r_{0})}$. Hence, $M$ has finite
topological type. Our proof is finished.
\end{proof}

\begin{corollary}
Theorem \ref{theorem2} is true.
\end{corollary}

 \end{document}